\documentclass[11pt]{article}

\usepackage[colorlinks=false,pdfborderstyle={/W 1}]{hyperref}

\def\arXiv#1#2{\href{http://arxiv.org/abs/#1}{{\tt arXiv:#1 [#2]}}}

\usepackage{amsmath}
\usepackage{amssymb}
\usepackage{amsthm}
\usepackage{amsfonts}
\usepackage{graphicx}
\usepackage{mathrsfs}

\usepackage[dvipsnames]{xcolor}

\newtheorem{theorem}{Theorem}[section]
\newtheorem{question}[theorem]{Question}

\newtheorem{conjecture}[theorem]{Conjecture}
\newtheorem{lemma}[theorem]{Lemma}
\newtheorem{proposition}[theorem]{Proposition}
\newtheorem{corollary}[theorem]{Corollary}
\theoremstyle{remark}
\theoremstyle{definition}
\newtheorem{definition}[theorem]{Definition}

\newcommand{\Z}{\mathbb{Z}}
\newcommand{\N}{\mathbb{N}}
\newcommand{\T}{\mathbb{T}}
\newcommand{\F}{\mathbb{F}}

\newcommand{\Cay}{\mathsf{Cay}}
\newcommand{\nh}{\mathsf{nh}}
\def\hent{\mathsf{h}}

\newcommand{\lora}{\longrightarrow}

\numberwithin{equation}{section}
\numberwithin{figure}{section}

\def\hitem#1#2{\item[\hypertarget{#1}{#2}]\expandafter\gdef\csname LBL#1ITM\endcsname{#2}}
\def\iref#1{\hyperlink{#1}{\csname LBL#1ITM\endcsname}}

\let\qqed=\qed
\def\QED{\qqed\medskip}
\let\qed=\QED

\def \eps {\varepsilon}

\def\P{\mathbb{P}}
\newcommand{\E}{\mathbb{E}}
\def\md{\mid}
\def\Bb#1#2{{\def\md{\bigm| }#1\bigl(#2\bigr)}}
\def\BB#1#2{{\def\md{\Bigm| }#1\Bigl(#2\Bigr)}}
\def\Bs#1#2{{\def\md{\mid}#1(#2)}}
\def\Pb{\Bb\P}
\def\Eb{\Bb\E}

\def\EB{\BB\E}
\def\Ps{\Bs\P}
\def\Es{\Bs\E}

\def\proofof#1{{ \medbreak \noindent {\bf Proof of #1.} }}

\def\bl{\begin{lemma}}
\def\el{\end{lemma}}
\def\bth{\begin{theorem}}
\def\eth{\end{theorem}}
\def\bc{\begin{corollary}}
\def\ec{\end{corollary}}
\def\bcj{\begin{conjecture}}
\def\ecj{\end{conjecture}}
\def\bpr{\begin{proposition}}
\def\epr{\end{proposition}}
\def\bde{\begin{definition}}
\def\ede{\end{definition}}
\newcommand{\be}{\begin{eqnarray}}
\newcommand{\ee}{\end{eqnarray}}
\newcommand{\bes}{\begin{eqnarray*}}
\newcommand{\ees}{\end{eqnarray*}}

\usepackage{mathrsfs}

\def\1{1\!\! 1}

\def\cA{\mathcal A}

\def\cC{\mathcal C}

\def\cR{\mathcal R}

\def\Unif{\mathsf{Unif}}
\def\Aut{\mathsf{Aut}}

\makeatletter 
\newcommand\mynobreakpar{\par\nobreak\@afterheading} 
\makeatother
\newenvironment{myitemize}{\mynobreakpar\begin{itemize}}{\end{itemize}}
\def\bit#1{
\ \vskip 4 pt 
\begin{myitemize}
\setlength{\parskip}{-2pt}
#1
\end{myitemize}
\vskip -7 pt
\nobreak
}

\title{Quantitative indistinguishability and\\
sparse and dense clusters in factor of IID percolations}
\author{
Endre Cs\'oka \and  P\'eter Mester \and G\'abor Pete
}
\date{\today}

\begin{document}
\maketitle

\begin{abstract} 
Chifan-Ioana (2010) implies that, for any factor of IID percolation on any nonamenable Cayley graph $G$, there is a countable set of (strong) indistinguishability classes for non-hyperfinite clusters. We introduce here quantitative strengthenings, called properties (qI) and (qSI): for clusters that are $\eta$-non-hyperfinite for some $\eta>0$, there are at most $M(G,\eta)<\infty$ (strong) indistinguishability classes, for any FIID percolation. 

We first show that (qI) and (qSI) for a Cayley graph $G$ are equivalent to the ``sparse implies thin'' property (SiT): any FIID percolation with $\eta$-non-hyperfinite clusters has density at least $c(G,\eta)>0$. Also, (SiT) is independent of the finite generating set of a group. We then prove, using entropy inequalities, that (SiT) holds for free groups, even for weak FIIDs. On the other hand, recent work of Jardón-Sánchez, Mellick, Poulin, and Wróbel implies that (SiT) fails for weak FIIDs on non-exact, i.e., not property (A) groups.

Furthermore, (SiT) implies that the Bernoulli graphing over any non-hyperfinite FIID cluster is strongly ergodic, and that indistinguishability for non-hyperfinite FIID clusters is equivalent to strong indistinguishability. These results actually follow from the work of Chifan-Ioana for {\it every} nonamenable Cayley graph, but with non-probabilistic proofs. 

We also prove, again using entropy inequalities, this time for {\it all} nonamenable Cayley graphs, that any FIID percolation with high enough expected degree must have a density close to 1, and there must be a single indistinguishability class of such clusters. On Kazhdan groups, there must be a single such cluster.

Our results have finite counterparts: in any large girth $d$-regular graph sequence, any FIID subgraph of average degree at least $2+\delta$ must have density at least $c(d,\delta)>0$. In the uniform random $d$-regular graph $G_{n,d}$, this holds for every subgraph of average degree at least $2+\delta$.

\end{abstract}

\section{Introduction and results}

Consider a Cayley graph $G=\Cay(\Gamma,S)$ of an infinite group $\Gamma$ with finite symmetric generating set $S$, 
 and a $\Gamma$-invariant {\bf site percolation} process on $G$: a random subset $\omega$ of $V(G)$, with law $\P$ invariant under the natural translation action of $\Gamma$. 
  The connected components of $\omega$ are called clusters, and the cluster of a vertex $x$ is denoted by $\cC_x$. A $\Gamma$-invariant percolation is called ergodic if every $\Gamma$-invariant event has probability 0 or 1. We will assume that the reader is familiar with the basics of group-invariant percolations; see \cite{BLPS, LPbook, PGG} for background. 

An ergodic invariant percolation process has {\bf indistinguishable infinite clusters}, defined by Lyons and Schramm \cite{LySch}, if, for every $\Gamma$-invariant Borel-measurable $\cA \subset \Omega := \{0,1\}^{V(G)}$, either  almost surely every infinite cluster of $\omega$ is in $\cA$, or almost surely none. Typical properties of interest, e.g., are transience of simple random walk on the cluster, or being one-ended, or the value of the critical density $p_c$ for having an infinite cluster in Bernoulli percolation, or being hit infinitely often by simple random walk on the Cayley graph itself; one only considers infinite clusters because finite clusters are obviously distinguishable by cardinality. An equivalent formulation is that the cluster $\cC_o$ of a fixed vertex $o$, conditioned on being infinite, together with the Cayley graph as a supergraph, is {\bf extremal as a unimodular random rooted graph} \cite{AldLy,Lovasz,PGG}. Or, in the language of orbit equivalence relations \cite{KechMill}: the {\bf cluster subrelation} on $\Omega\times\Omega$ corresponding to the invariant percolation (that is, $\omega \sim \gamma(\omega)$ if{f} $o$ and $\gamma(o)$ are connected in $\omega$), restricted to infinite clusters, is {\bf ergodic}: if $\cA \subset \Omega$ is an event that is a union of infinite equivalence classes, then either
$\cA$ or $\cA^c$ has probability zero.  See \cite{GabLy} or \cite{Seb} for these equivalences.


It was proved in \cite{LySch} that the infinite clusters in Bernoulli percolation on any Cayley graph, and more generally,  in any ergodic {\bf insertion tolerant} invariant percolation (meaning roughly that a site, conditionally on the entire configuration elsewhere, is still open with a positive probability) are indistinguishable. Besides the intrinsic interest, there are many applications (of the results or the techniques): uniqueneness monotonocity for Bernoulli percolation \cite{relentless}, understanding the uniqueness threshold \cite{Schon}, continuity of the percolation probability $\theta(p)$ for $p>p_c$ \cite{vdBK}, and the measurable cost of Kazhdan (T) groups \cite{HuPe}. 

However, there are many natural invariant percolation processes that do not have insertion tolerance, only some weaker tolerance under certain surgeries, but indistinguishability is still expected. Proved examples (which are bond percolations, $\omega \subset E(G)$, but all the definitions are analogous) are the Wired and Free Uniform Spanning Forests \cite{HuNa}, the Free Uniform and Minimal Spanning Forests assuming that they differ from their wired versions \cite{Timar}, and some processes with two-ended clusters, such as the interchange process and the $O(n)$ loop model \cite{EAPT}. Again there are several applications, both probabilistic \cite{HuPer} and ergodic theoretical \cite{GabLy,GM}.

Insertion tolerance and its relatives are natural generalizations to go beyond Bernoulli percolation. Another natural direction is to consider {\bf factor of i.i.d.~(FIID) processes}: those obtained as $\Gamma$-equivariant functions from $[0,1]^{V(G)}$ or $[0,1]^{E(G)}$ to $\{0,1\}^{V(G)}$ or $\{0,1\}^{E(G)}$. Here, full indistinguishability does not hold. For instance, consider the $4$-regular tree $\T_4$, assign to each vertex $v$ an independent label $U_v \sim \Unif[0,1]$, and declare it to be red if $U_v \in [0,1/2]$, blue if $U_v \in [3/5,1]$, uncoloured otherwise. Then consider the red and blue subtrees induced by the red and blue vertices, respectively. 
Using the fact that $p_c(\T_4)=1/3$  for Bernoulli percolation, it is easy to check \cite[Exercise 5.7]{LPbook} that the infinite red trees have $p_c=2/3$, while the blue ones have $p_c=5/6$. 
 A natural example where indistinguishability is not known (but expected) is the Wired Minimal Spanning Forest \cite{LPS,Timar}: there should be infinitely many trees in the forest for Cayley graphs of sufficiently high volume growth, but, if having an infinite cluster at $p_c$ has not been excluded (e.g., on most non-Abelian groups of subexponential growth, see \cite{HeHu} and the references there), then it is not known if all trees are 1-ended, or there are also some 2-ended ones.

An important theorem of Chifan and Ioana \cite{ChI} is that, for any equivalence subrelation $\cR$ of the orbit equivalence for the Bernoulli shift $\Omega=[0,1]^\Gamma$ over any countable nonamenable group $\Gamma$, the probability space decomposes as $\Omega = \bigcup_{i\ge 0} \Omega_i$, where each $\Omega_i$ is a measurable union of $\cR$-classes, the restriction of $\cR$ to $\Omega_0$ is {\bf hyperfinite} in the sense that it is an increasing union of equivalence relations with finite classes only, while the restriction to each $\Omega_i$, $i \ge 1$, is non-hyperfinite and ergodic, moreover, strongly ergodic: every asymptotically $\cR$-invariant sequence $\cA_n\subset \Omega_i$ of events satisfies $\Ps{\cA_n}\Ps{\cA_n^c}\to 0$. In particular, for the cluster subrelation of any FIID percolation on any Cayley graph, the theorem says that the non-hyperfinite clusters fall into countably many (strong) indistinguishability classes. In this setting, {hyperfiniteness} of $\cC_o$ has the probabilistic meaning that for any $\eta>0$ there is a jointly unimodular subset of the vertices of $\cC_o$ such that joint probability of $o$ being in the subset is at most $\eta$, but the complement in $\cC_0$ has finite clusters only.  See \cite[Theorem 8.5]{AldLy} for the equivalence with the previous definition. Furthermore, {\bf strong indistinguishability} is the translation of strong ergodicity to the cluster subrelation, made in~\cite{Seb}. Namely, we call $\cA_n(x)$, $n\ge 1$, $x\in V(G)$, an {\bf asymptotic cluster property} sequence if, for each $n$, it is a measurable subset of the possible clusters $\cC_x$, invariant under graph automorphisms of $\cC_x$ fixing $x$, which satisfies
\be\label{e.asymC}
\Pb{ \1_{\cA_n(x)} = \1_{\cA_n(y)} \md \cC_x = \cC_y \textrm{  is infinite} } \to 1 \quad \textrm{for any }x,y\in V(G). 
\ee
{Strong indistinguishability} then means that, for any asymptotic cluster property sequence,
\be\label{e.strongind}
\Pb{ \1_{\cA_n(x)} = \1_{\cA_n(y)} \md |\cC_x| = |\cC_y| = \infty }\to 1
\ee
also holds. A variation on the example in \cite{Seb} of an ergodic invariant bond percolation with indistinguishable but not strongly indistinguishable infinite clusters is the Wired Uniform Spanning Forest on the standard 4-regular Cayley graph of the free group $\F_2$, with edges labelled by the generators $a$ and $b$. Here each cluster is a one-ended tree, hence, if $\cC_x = \cC_y$, then the unique infinite self-avoiding paths, $\pi_x$ from $x$ and $\pi_y$ from $y$, coincide from some point on. Let $\cA_n(x)$ be the event that the majority of the first $2n+1$ edges on $\pi_x$ are labelled by the generator $a$ or $a^{-1}$. Generating the WUSF by Wilson's algorithm \cite{LPbook}, one can check that $\cA_n(x)$ is an asymptotic cluster property sequence~\eqref{e.asymC}, while vertices in different components will disagree with a uniformly positive property, hence~(\ref{e.strongind}) does not hold. It is not an accident that the clusters in this example are one-ended trees, hence hyperfinite: it follows from \cite[Proposition 6]{ChI} that, for FIID non-hyperfinite clusters, indistinguishability and strong indistinguishability are equivalent.

A countable upper bound on the number of indistinguishability classes requires non-hyper\-finite\-ness, as demonstrated in \cite{Mester} quite dramatically, where a FIID labelling of the vertices of the standard 4-regular Cayley graph of $\F_2$ was constructed with $\mathsf{Unif}[0,1]$ marginals, but vertices with the same label spanning infinite clusters --- by the Chifan-Ioana theorem, necessarily hyperfinite ones. (If one does not like distinguishability with the $\mathsf{Unif}[0,1]$ labels, one can modify the clusters using their labels, obtaining uncountably many indistinguishability classes for the clusters themselves.) One can also merge some of these clusters to get countably infinitely many distinguishable non-hyperfinite FIID clusters, showing that Chifan-Ioana's theorem is sharp also in this respect.

The proofs in \cite{ChI} use von Neumann algebras, building on rigidity results of \cite{Popa}. Earlier work of Ozawa \cite{Oza} uses $C^*$-algebraic techniques to prove the same result for exact groups, a notion we will define a few paragraphs later. Our initial motivation was twofold: 
\bit{
\item to provide a probabilistic approach to these indistinguishability results, at least in some special cases; 
\item to obtain new conditions that forbid an invariant percolation on a nonamenable Cayley graph to be a FIID; the currently known conditions are about correlation decay \cite{LyNaz,BaVi} and local entropy inequalities \cite{Bowen,Rahman,BSz,BGH,CsHaVi}. For instance, presently no tail-trivial invariant process on a regular tree is known that is not a FIID \cite{Russ}. Another, even more fantastic, application could be to show that on certain nonamenable Cayley graphs with {\bf measurable cost} 1, which means that for any $\delta>0$ there exists an invariant percolation with expected degree at most $2+\delta$ that makes the Cayley graph connected (very non-trivial examples are infinite Kazhdan groups \cite{HuPe}), such percolations might not be obtained in a FIID manner, hence they would not have {\bf fixed price} 1. This is a famous question of Gaboriau \cite{Gab}; see \cite{FMW,Ali} for recent examples with fixed price 1, and \cite{Zoo,Vor2,Vor5} for the related question of FIID Sparse Unique Clusters.
}
An unexpected turn of events was that the natural approach we have found actually addresses a strictly stronger property than the indistinguishability result of \cite{ChI}, as we will soon explain.

We call a unimodular random rooted graph (a URG), such as the cluster $\cC_o$ in an invariant site percolation $\omega$ on a Cayley graph,  {\bf $\eta$-non-hyperfinite} if every unimodular random subgraph of it that misses the root $o$ with a joint probability at most $\eta$ has infinite clusters. (Mind that the density $\eta$ is understood within the URG $\cC_o$, not within the Cayley graph; this is important when $\omega$ has a small density within the Cayley graph.) We will denote by $\nh(\cC_o)$ the supremum of values $\eta$ for which $\cC_o$ is $\eta$-non-hyperfinite. We are now ready to start stating our results.

\bth\label{t.equiv} {\bf (a)} 
For any finitely generated Cayley graph $G$, the following are equivalent:
\begin{itemize}
\hitem{qI}{\rm (qI)} {\bf ``Quantitative Indistinguishability''} In any FIID site percolation, for clusters that are $\eta$-non-hyperfinite for some $\eta>0$, there are at most $M(G,\eta)<\infty$ indistinguishability classes.
\hitem{qSI}{\rm (qSI)} {\bf ``Quantitative Strong Indistinguishability''} In any FIID site percolation, for clusters that are $\eta$-non-hyperfinite for some $\eta>0$, there are at most $M(G,\eta)<\infty$ strong indistinguishability classes.
\hitem{SiT}{\rm (SiT)} {\bf ``Sparse implies Thin''} Any FIID site percolation with only $\eta$-non-hyperfinite clusters has density at least $c(G,\eta)>0$. 
\end{itemize}
\vskip -0.2 cm
\bit{
\item[{\bf (b)}] 
If a group $\Gamma$ has a finitely generated Cayley graph $G=\Cay(\Gamma,S)$ that satisfies \iref{SiT}, then every finitely generated Cayley graph $H=\Cay(\Gamma,T)$ satisfies it.
\item[{\bf (c)}] Both parts (a) and (b) hold for weak limits of FIID processes (a.k.a.~weak FIIDs), as well.
}
\eth

This equivalence is interesting in both non-trivial directions. There may be operator algebraic techniques to prove~\iref{qI} for certain groups, while we will use probabilistic techniques to prove~\iref{SiT} on trees. (See Section~\ref{s.entropy}.)

\bth[Trees satisfy (SiT)]\label{t.TreeSiT}
{\bf (a)} If $G$ is a $d$-regular tree with $d\ge 3$, and $\Gamma$ is group acting transitively and unimodularly of $G$, then any weak limit of $\Gamma$-factor of IID site percolations satisfies~\iref{SiT}. {\bf (b)} In fact, \iref{SiT} holds for any typical process on the tree.
\eth

The most natural examples here for $\Gamma$ are $\Aut(\T_d)$; the free group $\F_r$ on $r\ge 2$ generators, with $d=2r$; or a $d$-wise free product $\Z_2^{\star d}$. In part (b), a {\bf typical process} is any Benjamini-Schramm limit of colourings on the uniform random $d$-regular graph $G_{n,d}$ on $n$ vertices \cite{BSz}. 

The proof of Theorem~\ref{t.equiv} is not difficult; the interest is more in the statement itself. The proof of Theorem~\ref{t.TreeSiT} uses {\bf entropy inequalities}, which originate from the works of Bollobás \cite{Boll}, Lewis Bowen \cite{Bowen} and Gamarnik-Sudan \cite{GaSu}, made explicit in \cite{BSz}. The original idea was that  FIID processes on regular trees can be approximated by processes on uniformly random finite regular graphs, which are locally tree-like, and where the local neighbourhood statistics of any process conform to some entropy inequalities that come from enumerating possible configurations in the finite graphs. A different approach was found in \cite{CsHaVi}, proving generalizations of these entropy inequalities for FIID processes on any nonamenable Cayley graph, using a proof similar to the construction of the Shapley value in cooperative game theory and Shearer's entropy inequality.

Beyond trees, how far could this result extend? It is not clear to us how to use the generalization in~\cite{CsHaVi} for this purpose. The big advantage of trees is that, for infinite unimodular random trees, $\eta$-non-hyperfiniteness is equivalent to having an expected degree $2+\delta$, with $\delta \asymp \eta$ (see \cite{BLPS}), and this degree condition, being completely local, is easy to use in entropy inequalities. The generator independence is a small step beyond trees, but, unfortunately, the constant $c(G,\eta)$ depends on the Cayley graph $G$, not on the group.

A counterexample by Tom Hutchcroft to a related but independent question of Miko{\l}aj Fr\c{a}czyk, at a workshop in Budapest in June 2024, made us realize that~\iref{SiT} may not hold for {\bf non-exact groups}. Exactness for groups, through their $C^*$-algebras, was introduced in \cite{exact}. Independently, property (A) for infinite graphs was defined in \cite{Yu}, where it was proved that these graphs have so-called uniform metric embeddings into Hilbert space. By \cite{ULA}, for finitely generated Cayley graphs, {\bf property (A)} is the existence of a connected subgraph $A_n(x)$ for each $n\in \N$, $x\in V(G)$, which contains $x$, such that $\big| A_n(x) \triangle A_n(y) \big| / \big( |A_n(x)| + |A_n(y)| \big) \to 0$ as $n\to\infty$, for every $x,y \in V(G)$. The equivalence of exactness and property (A) was noticed in~\cite{HiRo}; see \cite{AD} for an overview. The existence of groups without uniform embeddings into Hilbert space, because of an expander subgraph sequence, was first suggested in \cite{Gromov}, proved in \cite{Osajda}; these groups are nowadays called Gromov-Osajda monsters. In fact, not having property (A) for a finitely generated group has been shown to be equivalent \cite{ULA, Elek} to having {\bf small scale expanders} as induced subgraphs, meaning that there exists some $\kappa>0$ such that for every $n\ge 1$ there is a finite induced subgraph $G_n\subseteq G$ with the property that for every $S\subseteq V(G_n)$ with $|S|\leq n$ we have $|E(S,V(G_n)\setminus S)| > \kappa |S|$. These subgraphs were used recently in \cite{Krakow} to construct invariant random equivalence relations (IREs) on non-exact groups that are weak limits of finite IREs on the group, but are non-hyperfinite. It was also proved there that this is impossible on exact groups. Their construction easily implies the following:


\bth[Theorem 1.1 or 5.1 of \cite{Krakow} rephrased]\label{t.nonexact}
There are weak FIID site percolations on any countable non-exact group such that~\iref{SiT} does not hold.
\eth

Two obvious questions remain; we expect positive answers.

\begin{question}\label{q.exact} 
Does~\iref{SiT} hold for (weak) FIIDs on all finitely generated exact groups? 
\end{question}

\begin{question}\label{q.nonexact} 
Does~\iref{SiT} fail for some true FIIDs on non-exact groups? 
\end{question}

Given that Ozawa \cite{Oza} proved the Chifan-Ioana theorem only for exact groups, one may speculate that his method actually gives the stronger result~\iref{qI}. Very recently, Miko{\l}aj Fr\c{a}czyk has informed us that, with entropy arguments inspired by~\cite{CsHaVi} beyond trees, he has managed to give a positive answer to Question~\ref{q.exact} on all exact groups, even for weak FIIDs \cite{Mikolaj}.
\medskip

Despite Theorem~\ref{t.nonexact}, we can give an indistinguishability result that holds for (weak) FIIDs over every finitely generated group, again using the entropy inequality of \cite{CsHaVi}, see Section~\ref{s.entropy}.

\bth[Large degree implies high density]\label{t.lardegLV} 
For any $\eps,\gamma>0$ and integer $d\geq 3$, there exists $\delta=\delta(\eps,\gamma,d)>0$ such that, for any (weak) FIID site percolation $X$ on a Cayley graph with degree $d$ and external vertex Cheeger constant $\gamma$, if $\Es{\deg_X (o) \md X_o=1}\geq d(1-\delta)$ holds, then we also have that 
$\Ps{X_o=1}\geq 1-\eps$. 
\eth

This will immediately imply:

\bc[Large degree implies unique indistinguishability type]\label{c.lardegST} 
For any $\gamma>0$ and integer $d\geq 3$, there exists a $\delta=\delta(\gamma,d)>0$ such that in a Cayley graph with degree $d$ and external vertex Cheeger constant $\gamma$, for  any (weak) FIID site percolation $X$ on $G$, the infinite clusters whose average degree is at least $d(1-\delta)$ are indistinguishable.
\ec

A further corollary may be considered as a generalization of the result of \cite{LySch} that groups with Kazhdan's property (T) have a uniqueness threshold $p_u<1$ for Bernoulli percolation, and more generally, insertion-tolerant invariant percolations with high enough average degree have a unique infinite cluster. Even more generally, a high enough average degree of any invariant percolation on a Kazhdan group implies the non-vanishing of the two-point connectivity function, and in fact this property characterizes Kazhdanness \cite[Theorem 1.2]{MukhReck}.

\bc[High degree FIID clusters in Kazhdan groups are unique]\label{c.lardegKAZH} 
For any $d$-regular Cayley graph $G$ of a Kazhdan group, there exists a $\delta=\delta(G) > 0$, such that if $X$ is a (weak) FIID on $G$, then there can be at most one cluster whose average degree is at least $(1-\delta)d$. 
\ec

Let us stress that in a FIID percolation with average degree at least $d(1-\delta)$ on a Kazhdan group, besides the unique fat giant, there can also exist other, much thinner, infinite clusters. On the other hand, we have not found a counterexample to the following question:

\begin{question}[Fat giants on Kazhdan groups] Does the result of Corollary~\ref{c.lardegKAZH} hold for general invariant percolations on Kazhdan groups?
\end{question}

There are some implications of our notion of~\iref{SiT} regarding strong ergodicity. They actually follow from \cite{ChI} in a larger generality, for all nonamenable Cayley graphs, but with proofs that are abstract and convoluted for us, and possibly for most of the probability community. Moreover, our proofs give quantitative bounds, which we suspect to hold only for groups satisfying~\iref{SiT}.

\bpr[The Bernoulli graphing over (weak) FIID clusters is a large-scale expander]\label{pr.expander} If $G$ satisfies \iref{SiT}, and $\cC_o$ is an ergodic $\eta$-non-hyperfinite URG given by a (weak) FIID site percolation on $G$, then the Bernoulli graphing over $\cC_o$ is a large-scale expander (a.k.a.~strongly ergodic): for any FIID two-colouring of the vertices of $\cC_o$, with one-point marginal in $(\eps,1-\eps)$, the density of edges with endpoints of different colours is at least $\delta(G,\eta,\eps)>0$.
\epr

Here, the one-point marginal is understood in the annealed sense, as the probability of a colour in the joint probability space of $\cC_o$ and its FIID colouring. Then, by the joint ergodicity of $\cC_o$ and its FIID colouring, the density of a colour (measured, e.g., along a delayed simple random walk, stationary for $\cC_o$) is an almost sure constant in the probability space of $\cC_o$. (See the proof for more details.) A similar statement holds for the edge-density, by which we mean the probability that $o$ and a uniform random neighbour of its has different colours. 
A straightforward corollary is the following:

\bc[Indistinguishability implies strong indistinguishability]\label{c.strong}  
If $G$ satisfies~\iref{SiT}, and a (weak) FIID site percolation in it has indistinguishable clusters, then any asymptotic cluster property sequence $\cA_n(\cdot)$ is asymptotically trivial: $\Ps{ \cA_n(x) \,|\, |\cC_x|=\infty} \, \Ps{ \cA_n(x)^c \,|\, |\cC_x|=\infty } \to 0$. In particular, indistinguishable infinite clusters are strongly indistinguishable.
\ec

It was possibly first proved in \cite{almost,Klaus}, popularized by \cite{LyNaz}, that the Bernoulli graphing over any nonamenable Cayley graph $G$ is an expander (meaning that the mixed edge density is at least $\delta(G,\eps) \ge c\, \eps$, with a constant $c=c(G)>0$). Thus one may naively think that this holds over any non-hyperfinite URG as well, at least when the graph is almost surely nonamenable. This is not the case, as was pointed out by Miko{\l}aj Fr\c{a}czyk, presented in \cite{Halloween}: there exist almost surely nonamenable unimodular random graphs where the Bernoulli graphing is not even a large-scale expander. This shows that Proposition~\ref{pr.expander} cannot be entirely trivial. Moreover, the same example, together with the next proposition, imply that the entropy inequalities of~\cite{CsHaVi} cannot hold in the generality of almost surely nonamenable URGs.

\bpr[Entropy inequality and Bernoulli large-scale expansion]\label{pr.expander2}
If $(G,o)$ is an ergodic URG satisfying the entropy inequality~(\ref{e.starURG}), then its Bernoulli graphing is a large-scale expander.
\epr

Finally, Theorem~\ref{t.TreeSiT} has the following counterparts for finite graphs:

\bth[``Sparse implies Thin'' in $G_{n,d}$ and other large girth graphs]\label{t.finite}\ \\
{\bf (a)} In any essentially large girth $d$-regular graph sequence, any FIID subgraph of average degree at least $2+\delta$ must have density at least $c(d,\delta)>0$.\\
{\bf (b)} In the uniform random $d$-regular graph $G_{n,d}$ on $n$ vertices, any subgraph  of average degree at least  $2+\delta$ must have density at least $c(d,\delta)>0$.
\eth

Theorem~\ref{t.lardegLV} also has a finite counterpart: in any $d$-regular $\gamma$-expander sequence, or any $\eps>0$ there exists a $\delta = \delta(d,\gamma,\eps)>0$ such that
any subgraph of average degree at least $(1-\delta)d$ must have density at least $1-\eps$. This is actually very easy from the definition of an expander --- an exercise for the reader. And this implies Theorem~\ref{t.lardegLV} for Cayley graphs $G$ that have sofic (Benjamini-Schramm \cite{AldLy}) approximations by expanders $G_n$, since any (weak) FIID process on $G$ can be modelled with a small error on $G_n$. However, Corollary~\ref{c.lardegKAZH} is not implied even in the sofic case: a unique giant cluster on the finite graphs may fall apart into many infinite pieces in the limit.
\medskip

The rest of the article is organized as follows. Section~\ref{s.genproofs} contains the proofs of the general equivalences and their corollaries, while Section~\ref{s.entropy} contains the proofs that use entropy inequalities, namely those of the two key results, Theorems~\ref{t.TreeSiT} and~\ref{t.lardegLV}, and of their relatives, Proposition~\ref{pr.expander2} and Theorem~\ref{t.finite}. Note, in particular, that some corollaries of Theorems~\ref{t.TreeSiT} and~\ref{t.lardegLV} will hence be proved in the earlier Section~\ref{s.genproofs}.

\section{Proofs of the general equivalences and their corollaries}\label{s.genproofs}

\proofof{Theorem~\ref{t.equiv}} The implication \iref{qSI} $\implies$ \iref{qI} is trivial. The opposite direction follows immediately from Corollary~\ref{c.strong}.
\medskip

Assume that \iref{qI} does not hold: for any $M>0$, there are more than $M$ indistinguishability classes for $\eta$-non-hyperfinite clusters (possibly infinitely many, even uncountably). If there exist $M$ such classes of positive probability each, then at least one of them has probability in $(0,1/M]$. Otherwise, there is some positive probability taken by non-atomic classes. In either case, there exists some $\Gamma$-invariant event of probability in $(0,1/M]$ that an $\eta$-non-hyperfinite $\cC_o$ can satisfy; keeping all infinite clusters with this property and deleting everything else, we get a FIID percolation with only $\eta$-non-hyperfinite clusters and density at most $1/M$. If $M >  1/c(G,\eta)$, we get a contradiction to~\iref{SiT}.
\medskip

Now assume that \iref{SiT} does not hold; i.e., there are FIID percolations with $\eta$-non-hyperfinite clusters and densities $\eps_n \searrow 0$. By passing to a subsequence, we may assume that $10 d \sum_{n\ge 1} \eps_n < \eta$, where $d$ is the degree of $G$. Take independent copies $\omega_n \subset V(G)$ of each FIID percolation. Let $\overline\omega_n$ be the union of $\omega_n$ with its exterior vertex boundary in $G$. Since $\omega_n$ is independent of $\bigcup_{i\not=n} \overline\omega_i$, the density of $\omega_n \cap \bigcup_{i\not=n} \overline\omega_i$ inside $\omega_n$ is less than $\eta/2$, and hence, by $\omega_n$ being $\eta$-non-hyperfinite, $\omega_n \setminus \bigcup_{i\not=n} \overline\omega_i$ has some infinite clusters $\cC$, with $p_c(\cC) < 1-\eta/2$. Let $\omega^*_n$ be the union of all these infinite clusters $\cC$. For different $n$'s, these $\omega^*_n$'s are all disjoint. Pick an infinite sequence of distinct reals $r_1,r_2,\dots \in (1-\eta/2,1)$. By applying Bernoulli percolation at density $p_c(\omega^*_n)/r_n < 1$ to $\omega^*_n$, we get some infinite clusters with $p_c=r_n$. The union of these clusters is a FIID percolation with infinitely many disitinguishable infinite clusters, hence \iref{qI} does not hold. 
\medskip

For the proof of {\bf (b)}, assume that $\omega$ is a $\Gamma$-invariant FIID site percolation on $H$ with $\eta$-non-hyperfinite clusters only. Let $\xi$ be the following site percolation on $G$: if $x,y$ are neighbours in $H$ and are both open in $\omega$, then take all the finitely many shortest paths in $G$ between $x$ and $y$, and open all of their vertices. This is clearly a $\Gamma$-invariant FIID process, because for any $v \in E(G)$, the set $E_{G,H}(v)$ of $(x,y)\in E(H)$ for which $v$ lies on a shortest $x$-to-$y$ path in $G$ is a $\Gamma$-equivariant finite set.

We claim that  every cluster of $\xi$ is $\eta/K$-non-hyperfinite, where $K=K(\Gamma,S,T)=|E_{G,H}(v)|<\infty$ from the previous paragraph. Indeed, first look at the vertices of $\cC_o(\omega)$, originally a subgraph of $H$, now inside $G$, together with all the shortest paths between them, which are all open in $\xi$ by definition. Denote this subgraph of $G$ by $\overline{\cC_o(\omega)}^{\xi}$; this might be a strict subset of the full cluster $\cC_o(\xi)$, since some of these $\overline{\cC_u(\omega)}^{\xi}$, with $u\not\in \cC_o(\omega)$, may intersect $\overline{\cC_o(\omega)}^{\xi}$. This $\overline{\cC_o(\omega)}^{\xi}$ is a unimodular random subgraph of $G$. Assume that there exists a unimodular subset of its vertices with marginals at most $\eps>0$ that breaks it into finite components. Each vertex $v$ in $\overline{\cC_o(\omega)}^{\xi}$ has at most $K$ ``preimages'' in $E_{G,H}(v) \cap \cC_o(\omega)$; if we delete all of these from $\omega$, then $\cC_o(\omega)$ will fall apart into finite components. Thus $\eps K \ge \eta$, meaning that  $\overline{\cC_o(\omega)}^{\xi}$ is $\eta/K$-non-hyperfinite. And then, $\cC_o(\xi)$ being a unimodular supergraph of $\overline{\cC_o(\omega)}^{\xi}$, it is also $\eta/K$-non-hyperfinite.

By $G$ satisfying \iref{SiT}, the density of $\xi$ in $G$, hence in $H$, is at least $c(G,\eta/K)>0$, thus we have \iref{SiT} also for $H$.
\medskip

Part~{\bf (c)} is clear, since all the extra constructions we used in the proofs are themselves FIIDs.
\qed

\proofof{Theorem~\ref{t.nonexact}}
For any Cayley graph $G$ of a countable nonexact group, recall from the Introduction the sequence of finite induced subgraphs $G_n\subset G$, the small scale expander sequence for some $\kappa>0$. The basic construction of \cite[Thm 5.1]{Krakow} is that for each $v\in V(G)$ we place a translated copy of $G_n$ with probability $\eps/|V(G_n)|$, independently. The density of this FIID site percolation is $\asymp \eps$, see \cite[Lemma 5.3]{Krakow}. Now take any weak limit point of these percolations; that will have density $\asymp \eps$, and will be $\eta(\kappa,d)$-non-hyperfinite by the same combination of Propositions 5.2 and 5.4 of \cite{Krakow} that gave their Theorem 5.1.
\qed

\proofof{Corollary~\ref{c.lardegST}} Choose $\eps<1/2$, and let $\delta$ be the corresponding constant guaranteed by Theorem~\ref{t.lardegLV}. If there was a property $\cA$ that holds for some but not all clusters of degree at least $d(1-\delta)$, then we get a contradiction, as follows. The condition of Theorem~\ref{t.lardegLV} would still hold for the following two percolations: $X^{\tt yes}$ which only keeps the clusters whose average degree is at least $d(1-\delta)$ and for which $\cA$ holds, and $X^{\tt no}$
which only keeps the clusters whose average degree is at least $d(1-\delta)$ and for which $\cA$ does not hold. The vertex sets of these clusters are disjoint, but by the theorem they both must have density at least $1-\eps$, which is impossible.\qed

\proofof{Corollary~\ref{c.lardegKAZH}} Recall that a group has Kazhdan's property (T) if{f} for any of its finitely generated Cayley graphs $G$ and $\eps>0$ there is some $\kappa(G,\eps)>0$ such that any ergodic two-colouring of the vertices with marginals in $(\eps,1-\eps)$ has at least a density $\kappa$ of mixed coloured edges \cite{GlaW}. 

Let us call a cluster {\it $\delta$-fat} if the average degree within it is at least $d(1-\delta)$.
By the ergodicity of $X$, the number of $\delta$-fat clusters is an almost sure constant $K_{\delta}\in \N\cup\{\infty\}$.
If $K_{\delta}=\infty$, then colour every $\delta$-fat cluster with fair bits from $\{0,1\}$, independently for the different $\delta$-fat clusters, and let a vertex belonging to a $\delta$-fat cluster inherit the bit of its cluster (so the colouring is constant on the $\delta$-fat clusters). The vertices which are not included in one of the $\delta$-fat clusters can be coloured independently by fair bits. This colouring is ergodic, because the infinitely many $\delta$-fat clusters are all indistinguishable by Corollary~\ref{c.lardegST}, hence each has frequency zero for simple random walk on $G$, and hence \cite[Section 2.3]{HuPe} applies. The marginal of each colour is $1/2$, so the density of mixed edges must be at least $\kappa(1/2)$. On the other hand, as $\delta\to 0$, Theorem~\ref{t.lardegLV} easily implies that the density of mixed edges also goes to $0$, which is a contradiction. If $K_{\delta}=k>1$, then pick $\lfloor k/2 \rfloor$ of the $\delta$-fat clusters uniformly at random and colour them by $1$, while colour the remaining ones with $0$. As before, colour the vertices not belonging to any of these clusters with i.i.d.~fair bits again. By some strange coincidence, it may happen that this colouring is not ergodic --- in that case, take its ergodic decomposition. In each ergodic component, the marginal of the colours will be in $[\frac{1}{3},\frac{2}{3}]$. So, the density of mixed edges must be at least $\kappa(1/3)$, but as before, as $\delta\to 0$, the density of mixed edges must go to $0$ as well. So, we have arrived at a contradiction unless $K_{\delta}\in\{0,1\}$.
\qed

\proofof{Proposition~\ref{pr.expander}} Assume that, for the ergodic $\eta$-non-hyperfinite (weak) FIID cluster $\cC_o$, there is a sequence of FIID two-colourings $\xi_n: V(\cC_o) \lora \{0,1\}$ with marginals in $(\eps,1-\eps)$ but the density of mixed colour edges $\delta_n\to 0$. (Also for non-transitive graphs $G$, by a FIID colouring we mean a measurable map from IID labels on $V(G)$ or $E(G)$ that is $\Aut(G)$-equivariant.) Note that the probability that $o$ has a neighbour $v$ with $\xi_n(u)\not=\xi_n(v)$ is at most $d\delta_n$, if $G$ has degree $d$. By passing to a subsequence, we can assume that $\sum_{n\ge 1} d \delta_n  < \eta/2$. Then, if we remove every vertex $u$ from $\cC_o$ that has some neighbour $v$ in $\cC_o$ and some $n$ with $\xi_n(u)\not=\xi_n(v)$, then on each of the resulting components all the $\xi_n$s will be constant, and the resulting component $\cC^*_o$ of the origin will have $\Eb{\nh(\cC^*_o)} > \nh(\cC_o)-\eta/2>\nh(\cC_o)/2$. 

Consider now $c(G,\eta/2)>0$ from \iref{SiT}. If we can ensure that, for some large enough $\ell$ and well-chosen colourings $\xi_1,\dots,\xi_\ell$, the density of sub-clusters for each of the $2^\ell$ colour-classes is smaller than this $c(G,\eta/2)$, then we arrive at a contradiction, since each colour-class defines a (weak) FIID percolation, and at least one of them must have $\nh(\cC^*_o) > \eta/2$ (since the expectation over all colour classes is larger than $\eta/2$). 

So, in the rest of the proof, we will choose a further subsequence of the colourings; by a slight abuse of notation, they still will be denoted by $\xi_1,\xi_2,\dots$. 




Consider the delayed simple random walk $(X_n)_{n\ge 0}$. For any fixed $k\in\N$, the sequence 
$$(\cC_o,X_n,\xi_1(X_n),\xi_2(X_n),\dots,\xi_k(X_n))_{n\ge 0}$$ is stationary in $n$, moreover, ergodic. See~\cite[Lemma 2.2]{R3T3} for a proof. (And note that it is not important here that $\cC_o$ is itself a FIID.)




With $\Ps{\xi_1(o)=1}=\eps_1 \in (\eps,1-\eps)$, by Birkhoff's ergodic theorem, $\frac{1}{N} \sum_{n=0}^{N-1} \xi_1(X_n) \to \eps_1$ almost surely, independently of $\xi_1(o)$. Hence there must exist some deterministic $n$ with 
\be\label{e.Birkhoff}
\Pb{ \xi_1(X_0)=1, \  \xi_1(X_n)=0 } \geq 0.99 \,  \eps_1(1-\eps_1) > \eps/4.
\ee

Fix such an $n$. There exists a $\xi_2$ such that its mixed edge density $\delta_2$ satisfies $n \delta_2 < \eps/4$. Then, 
$$2\, \Pb{ \xi_2(X_0)=1, \  \xi_2(X_n)=0 } = \Pb{ \xi_2(X_0) \not= \xi_2(X_n) } \leq \sum_{i=1}^n \Pb{  \xi_2(X_{i-1}) \not= \xi_2(X_{i}) } \leq n\delta_2 < \eps/4.$$

By stationarity, $\Pb{\xi_1(o)=\xi_2(o)} = \Pb{\xi_1(X_0)=\xi_2(X_0)} = \Pb{\xi_1(X_n)=\xi_2(X_n)}$. Therefore,
\begin{align*}
\Pb{ \xi_1(X_0)=1,\ \xi_1(X_n)=0 } &\leq  \Pb{ \xi_1(X_0)=1,\  \xi_2(X_0)=0 } \\
& \hskip 1 cm + \Pb{ \xi_2(X_0)=1,\  \xi_2(X_n)=0 } + \Pb{ \xi_2(X_n)=1,\  \xi_1(X_n)=0 }\\
& \leq  \Pb{\xi_1(o) \not= \xi_2(o)} + \eps/8\,.
\end{align*}
Combining with \eqref{e.Birkhoff}, we get that 
$$
\eps/8 < \Pb{\xi_1(o) \not= \xi_2(o)}.
$$

By considering $\xi_2':=1-\xi_2$, we also obtain
$$
\eps/8 < \Pb{\xi_1(o) = \xi_2(o)}.
$$

These lower bounds, combined with every one-point marginal being in $(\eps,1-\eps)$, easily give that out of the four possible values of $\xi_{1,2}:=(\xi_1,\xi_2)$, at least three has density at least $\eps/8$ each. 

We can now repeat the previous argument with $\xi_{1,2}$ and a suitably chosen $\xi_3$. Namely, Birkhoff now gives
$$
\Pb{ \xi_{1,2}(X_0)=\mathsf{x}, \  \xi_{1,2}(X_n)=\mathsf{y} } \geq 0.99 \,  (\eps/8)^2 > \eps^2/100,
$$
for any $\mathsf{x}, \mathsf{y} \in \{0,1\}^2$ that have density at least $\eps/8$, with some $n=n_{\mathsf{x},\mathsf{y}}$. We then choose $\xi_3$ with $n \delta_3 < \eps^2/100$, and using that 
\begin{align*}
\Pb{ \xi_{1,2}(X_0)=\mathsf{x},\ \xi_{1,2}(X_n)=\mathsf{y} } &\leq  \Pb{ \xi_{1,2}(X_0)=\mathsf{x},\  \xi_3(X_0)=0 } \\
& \hskip 0.6 cm + \Pb{ \xi_3(X_0)=1,\  \xi_3(X_n)=0 } + \Pb{ \xi_3(X_n)=1,\  \xi_{1,2}(X_n)=\mathsf{y}}\,,
\end{align*}
and that second term on the RHS is at most $\eps^2/200$, we obtain
$$
\eps^2/200 < \Pb{\xi_{1,2}(o)=\mathsf{x},\ \xi_3(o)=0} + \Pb{\xi_{1,2}(o)=\mathsf{y},\ \xi_3(o)=1}.
$$
This is for any $\mathsf{x}, \mathsf{y} \in \{0,1\}^2$ that have density at least $\eps/8$, of which there are at least 3, as we showed above; call them $\mathsf{x}_1, \dots, \mathsf{x}_{k_2}$ with $k_2 \ge 3$. Just as before, our lower bounds easily imply that, for any such $\mathsf{x}\not=\mathsf{y}$, at least three of the four possibilities $(\mathsf{x},i)$, $(\mathsf{y},i)$, $i=1,2$, have probabilities at least $\eps^2/200$. Applying this to all the four-tuples formed by $(\mathsf{x_j},i)$, $j=1,\dots,k_2$, $i=1,2$, we get that at most one of the $2k_2$ pairs $(\mathsf{x_j},i)$ can fail to have probability at least $\eps^2/200$. So, out of the 8 possibilities for $(\xi_{1,2},\xi_3) \in \{0,1\}^3$, there are $k_3 \ge 2k_2-1 \ge 5$ that has density at least $\eps^2/200$ each; the others have smaller densities.

In the next step, we get that, out of the 16 possibilities in $\{0,1\}^4$, at least $k_4 \ge 2k_3-1 \ge 9$ has density at least $0.99 \, (\eps^2/200)^2 > \eps^4/50000$ each. Continuing in this manner, within the $2^\ell$ possibilities for $(\xi_1,\dots,\xi_\ell)$, there will be at least $k_\ell \ge 2k_{\ell-1}-1 \ge 2^{\ell-1}+1$ colour classes with density at least some $\eps_\ell > 0$ that converges to 0 quickly. But this means that each of these ``large'' densities is at most $1/k_\ell < 1/2^{\ell-1}$. So, the density of every single colour class of the $2^\ell$ possibilities is at most $\max\{1/2^{\ell-1},\eps_\ell\}\to 0$. For $\ell$ large enough, this will be smaller than $c(G,\eta/2)>0$ from the beginning of the proof, and we are done.\qed

\proofof{Corollary~\ref{c.strong}} Assume that there is a (weak) FIID site percolation with indistinguishable infinite clusters, and there is a nontrivial asymptotic cluster property sequence $\cA_n(x)$; by taking a subsequence, we can assume that $\lim_{n\to\infty} \Ps{\cA_n(o) \,|\, |\cC_o|=\infty}=\eps \in (0,1)$. Consider the two-colouring $\xi_n$ of the vertices $v\in\cC_o$ given by $\1_{\cA_n(v)}$. This a FIID (with no extra randomness beyond $\cC_o$), where $\Ps{\xi_n(o)\not=\xi_n(v)}\to 0$ for any neighbour $v$ of $o$, by $\cA_n(\cdot)$ being an asymptotic cluster property. However, this contradicts Proposition~\ref{pr.expander}.
\qed

\section{Using entropy inequalities}\label{s.entropy}

The proofs of Theorems~\ref{t.TreeSiT},~\ref{t.lardegLV},~\ref{t.finite} and Proposition~\ref{pr.expander2} all use the general entropy inequalities of \cite{CsHaVi}, for which the setting is the following.

Let $(X_v)_{v\in V}$ be a FIID vertex process, with an arbitrary source and a countable output for the factor map. Let $H(X_U)$ denote the $\log_2$-entropy. (All inequalities will be linear in $H$, hence are independent of the choice of the base.) We will also use the conditional entropy $H(X \,|\, Y)= H(X,Y)-H(Y)$, which is also the expected entropy of $X$, given $Y$. Let $U\mapsto \alpha_U \ge 0$, for finite sets $U\subset V(G)$, be a $\Gamma$-invariant weight function. Then \cite[Theorem 3.3]{CsHaVi} says that 
\be\label{e.genbeta}
\sum_{U \ni o} \frac{\alpha_U}{|U|} H(X_U) \ge \beta H(X_o)\,, \quad \text{with}\quad \beta := \inf_{W\subset V}\frac{\sum_{U: U\cap W\not=\emptyset} \alpha_U}{|W|}\,.
\ee
We will use this when $\alpha_U=1$ if $U$ is a star (including its center), $0$ otherwise. Then the LHS is $d+1$ times $\frac{H(X_S)}{d+1}$, where $S$ is the star of $o$. On the RHS,  $\beta=\beta_S = \inf |\overline W|/|W|$, where, as before, $\overline{W}:=W \cup \partial_V^{\rm{ext}} W$. Therefore,
\be\label{e.starbeta}
H(X_S) \ge  \inf_{W\subset V} \frac{|\overline W|}{|W|} \cdot H(X_o)\,.
\ee

On any nonamenable Cayley graph, $\beta_S =\inf |\overline W|/|W| = 1+\gamma$, with $\gamma>0$ being the {\bf external vertex Cheeger constant}.
So in general we have:
\be\label{e.starGEN}
H(X_S) \ge   (1+\gamma) \cdot H(X_o)\, ,
\ee
while, specializing this to the $d$-regular tree, it is easy to see that $\beta_S=d-1$, hence
\be\label{e.starvert}
H(X_S) \ge (d-1) H(X_o)\,.
\ee

Since these entropy inequalities are local, they survive taking weak limits, so also hold for weak FIID processes.

\proofof{Theorem~\ref{t.TreeSiT}~(a)} Now assume that we are on $\T_d$, and $X_o\in \{0,1\}$, and 
$$\Ps{X_o=1}=\eps \quad\text{and}\quad \Es{\deg (o) \md X_o=1}=2+\delta$$
 for some $\eps,\delta>0$. Then 
 $$
 H(X_o) = \hent(\eps):= \eps\log\frac1\eps + (1-\eps)\log \frac{1}{1-\eps}  \sim \eps\log(1/\eps)
 $$ 
 as $\eps\to 0$ (with the usual meaning of $\sim$ that the ratio of the two sides tends to 1), which gives an estimate for the RHS of~\eqref{e.starvert}.
To estimate the LHS of~\eqref{e.starvert}, let us write 
\be\label{e.star}
 H(X_S) = H(X_o)+ H(X_{S\setminus\{o\}} \,|\, X_o), 
\ee
and {express} the second term {using} the two possibilities:
\be\label{e.emptystar}
\begin{aligned}
H(X_{S\setminus\{o\}} \,|\, X_o) &= \eps \cdot H(X_{S\setminus\{o\}} \,|\, X_o=1) + (1-\eps) \cdot H(X_{S\setminus\{o\}} \,|\, X_o=0) \\
 \end{aligned}
\ee
We have the estimate $H(X_{S\setminus\{o\}} \,|\, X_o=1)\leq d$, simply because that is the maximal value the entropy can take when the image set is of size at most $2^d$.
If we introduce
\be\label{e.kappaINTRO}
\kappa:= \Es{X_\sigma \md X_o=0}, 
\ee
for $\sigma$ chosen uniformly at random from $S\setminus\{o\}$,
then we can estimate 
\be\label{e.sigmaent}
H(X_{S\setminus\{o\}} \,|\, X_o=0) \leq d \cdot H(X_\sigma \,|\, X_o=0)= d \, \hent(\kappa),
\ee
for instance because we can reveal the values in $X_{S\setminus\{o\}}$ one-by-one in a uniform random order, and in each step the additional conditional entropy is at most $H(X_\sigma \,|\, X_o=0)$. Luckily, we can find $\kappa$ using the assumptions:
\begin{align*}
d\,\eps = \EB{\sum_{v\in S\setminus\{o\}} X_v} &= \eps \cdot \Eb{\deg(o) \md X_o=1} + (1-\eps) \cdot \Eb{\deg(o) \md X_o=0} \\
&= \eps \, (2+\delta) + (1-\eps) \, d \, \kappa,
\end{align*}
thus we get
\be\label{e.kappa}
\kappa=  \frac{(d-2-\delta)\eps}{d(1-\eps)}\,.
\ee 


Combining (\ref{e.starvert}), (\ref{e.star}), (\ref{e.emptystar}), (\ref{e.sigmaent}) and~(\ref{e.kappa}) into one chain, we get the following: $$(d-1)H(X_o)\leq H(X_S)\leq H(X_o)+ H(X_{S\setminus\{o\}} \,|\, X_o)\leq H(X_o)+ \eps \cdot d + 1\cdot d \cdot \hent(\kappa).$$
Using that $\kappa\to 0$ as $\eps\to 0$, hence $\hent(\kappa)\sim \kappa \log(1/\kappa)$, we get
$$
(d-2) H(X_o) \leq \eps\, d + (d-2-\delta) (1+o_{\eps\to 0}(1)) \,  \eps \log (1/\eps) + O_d(\eps)\,.
$$
Rearranging,
$$
\delta\, \eps \log (1/\eps) \leq O_d(\eps)\,.
$$
As we let $\eps\to 0$, we see that $\delta\to 0$ must happen, proving~\iref{SiT}.

For part~{\bf (b)}, one should use the same entropy inequality~\eqref{e.starvert} proved for typical processes in \cite[Theorems 3 and 4]{BSz}.
\qed

\proofof{Theorem~\ref{t.lardegLV}} This time, we assume that 
$$\Ps{X_o=1}=1-\eps \quad\text{and}\quad \Es{\deg (o) \md X_o=1}=d(1-\delta)$$
 for some $\eps,\delta>0$. Similarly to the above, but slightly differently, we now estimate 
\be\label{e.emptystarHigh}
\begin{aligned}
H(X_{S\setminus\{o\}} \,|\, X_o) &= (1-\eps) \cdot H(X_{S\setminus\{o\}} \,|\, X_o=1) + \eps \cdot H(X_{S\setminus\{o\}} \,|\, X_o=0) \\
&\leq 1\cdot d \cdot \hent(\delta) + \eps \cdot d \cdot \hent(\kappa),
 \end{aligned}
\ee
where, for $\sigma$ chosen uniformly at random from $S\setminus\{o\}$, 
\begin{align}
\delta &= \Es{1-X_\sigma \md X_o=1}\,, \nonumber\\ 
\kappa &:= \Es{1-X_\sigma \md X_o=0} = 1+ \delta - \frac{\delta}{\eps} \,, \label{e.kappaHigh}
\end{align}
by a calculation as before:
\begin{align*}
d\,\eps = \E{\sum_{v\in S\setminus\{o\}} (1-X_v)} &= (1-\eps) \cdot \Eb{ d-\deg(o) \md X_o=1} + \eps \cdot \Eb{d - \deg(o) \md X_o=0} \\
&= (1-\eps)  \, d \, \delta + \eps \, d \,\kappa .
\end{align*}

Now take $\delta\to 0$. Then~\eqref{e.kappaHigh} gives
\be\label{e.hkappa}
\hent(\kappa) \sim (1/\eps-1) \, \hent(\delta) \sim (1/\eps-1) \, \delta\log(1/\delta)\,.
\ee
 

Thus, from~\eqref{e.starGEN}, \eqref{e.star}, \eqref{e.emptystarHigh} and~\eqref{e.hkappa}, we get
$$
\gamma \, \hent(\eps)  \leq H(X_{S\setminus\{o\}} \,|\, X_o) \leq  d \big( 2-\eps + o_{\delta\to 0}(1) \big) \hent(\delta)\,.
$$
It follows that, for all small enough $\delta>0$, we must have $\eps < \frac{2d}{\gamma} \delta \to 0$, finishing the proof.
\qed

\proofof{Proposition~\ref{pr.expander2}} The version of the star-vertex entropy inequality~\eqref{e.starGEN} that could make sense for almost surely $\gamma$-nonamenable URGs is this, we think:
\be\label{e.starURG}
H\big(X_S \,|\, (G,o)\big) \ge  (1+\gamma) \cdot H\big(X_o \,|\, (G,o)\big)\,.
\ee
With this, and with $d:=\Es{\deg_G(o)}$, the proof of Theorem~\ref{t.lardegLV} goes through verbatim, and we get that high degree FIID subgraphs have a high density.

Let us now assume that there exists a FIID two-colouring of $(G,o)$ with marginals in $(\eps,1-\eps)$ but mixed-edge density $\delta$ so small that the previous result implies that any FIID subgraph with average degree at least $d(1-\delta)$ has density larger than $1-\eps$. Then, as in the proof of Proposition~\ref{pr.expander}, removing all vertices that have a neighbour with a different colour, we arrive at two disjoint FIID site percolations, each with density less than $1-\eps$, a contradiction.\qed

\proofof{Theorem~\ref{t.finite}~(a)} Using the finite hypergraph version \cite[Theorem 2.4]{CsHaVi} in place of \cite[Theorem 3.3]{CsHaVi}, the proof is identical to that of Theorem~\ref{t.TreeSiT}. 

For part~{\bf (b)}, recall that the entropy inequality~\eqref{e.starvert} for typical processes in \cite{BSz} actually comes from counting arbitrary colourings in $G_{n,d}$, hence it holds also in this setting.
\qed

\section{Acknowledgements}

We thank Miklós Abért, Miko{\l}aj Fr\c{a}czyk, Viktor Harangi, Ben Hayes, Tom Hutchcroft, Héctor Jardón-Sánchez, Antoine Poulin, Ádám Timár, László Márton Tóth and Robin Tucker-Drob for conversations.

Our research was supported by the Hungarian National Research, Development and Innovation Office, grants KKP 138270 (Csóka) and Advanced 151155 (Pete); by the Icelandic Research Fund grant No.~239736-051 (Mester); and by the ERC grant No.~810115-DYNASNET (Pete).

\bibliographystyle{alpha}


\ \\
{\bf Endre Csóka}\\
HUN-REN Alfr\'ed R\'enyi Institute of Mathematics, Re\'altanoda u. 13-15, Budapest 1053 Hungary\\
\texttt{endre.csoka[at]renyi.hu}\\ \url{https://www.renyi.hu/en/staff/endre-csoka}\\
\ \\
{\bf Péter Mester}\\
Division of Mathematics, The Science Institute, University of Iceland\\
Dunhaga 3 IS-107 Reykjavik, Iceland\\
\texttt{pmester[at]alumni.iu.edu}\\
\ \\
{\bf G\'abor Pete}\\
HUN-REN Alfr\'ed R\'enyi Institute of Mathematics, Re\'altanoda u. 13-15, Budapest 1053 Hungary, and\\
Department of Stochastics, Institute of Mathematics, Budapest University of Technology and Economics, M\H{u}egyetem rkp.~3., Budapest 1111 Hungary\\
\texttt{gabor.pete[at]renyi.hu}\\
\url{http://www.math.bme.hu/~gabor}

\end{document}